\documentstyle[11pt]{article}

\begin{document}
\noindent
\begin{center}
  {\LARGE Surgery, quantum cohomology and birational geometry}
  \end{center}

  \noindent
  \begin{center}

    {\large Yongbin Ruan}\footnote{partially supported by a NSF grant and a
    Sloan
    fellowship}\\[5pt]
      Department of Mathematics, University of Wisconsin-Madison\\
        Madison, WI 53706\\[5pt]
(To appear in the Berkeley-Stanford-Santa Cruz symplectic geometry
conference proceedings)\\[5pt]

              \end{center}

              \def \J{{\cal J}}
              \def \Map{Map(S^2, V)}
              \def \M{{\cal M}}
              \def \A{{\cal A}}
              \def \B{{\cal B}}
              \def \C{{\bf C}}
              \def \Z{{\bf Z}}
              \def \R{{\bf R}}
              \def \P{{\bf P}}
              \def \I{{\bf I}}
              \def \N{{\cal N}}
              \def \T{{\cal T}}
              \def \Q{{\bf Q}}
              \def \D{{\cal D}}
              \def \H{{\cal H}}
              \def \S{{\cal S}}
              \def \e{{\bf E}}
              \def \CP{{\bf CP}}
              \def \U{{\cal U}}
              \def \E{{\cal E}}
              \def \F{{\cal F}}
              \def \L{{\cal L}}
              \def \K{{\cal K}}
              \def \G{{\cal G}}
              \def \k{{\bf k}}
\section{Introduction}
     Recently, some amazing relations between quantum cohomology and birational
geometry have been discovered. Surgeries play a fundamental role in these recent discoveries.
This expository article surveys these new relations.

In the early days of quantum cohomology theory, the author
\cite{R3} showed that Mori's extremal ray theory can be generalized to
symplectic manifolds from projective manifolds.
This gave the first evidence of the existence of a link between quantum cohomology
and birational geometry.
Subsequently, attention was diverted to
establishing a mathematical foundation for quantum cohomology.
The progress on this link was slow.
There are several papers (\cite{Wi2},\cite{Wi3} and others)
extending
results of \cite{R3} to Calabi-Yau 3-folds and crepant resolutions.
It is time to think these mysterious connections again.
The breakthrough appeared in  the paper of Li-Ruan \cite{LR} where they  calculated
the change of quantum cohomology under flops and small extremal transitions for 3-folds.
Their results showed a surprising
naturality property of quantum cohomology under these surgeries.
At present, this naturality property
is best understood in complex dimension three.
The scope of this new naturality property
is still unknown.
The author believes that there is some very interesting mathematics
hidden in these intrigue properties and that the general machinery of quantum cohomology can be very
useful in uncovering them.

One of the goals of this paper is to attract more people to work on this topic.
Hence, the author shall try to make some conjectures and proposals
for general cases to entertain the readers.
Since it is a survey instead of a research paper,  the author will concentrate on the ideas
and techniques. The main reference of this article is Li-Ruan's paper \cite{LR}.
The reader can find more details and more complete references there. The author apologizes for the lack of details.

    This paper is organized as follows. We give a brief survey on
quantum cohomology and the naturality problem in section~2 and fix  notation.
Then, we will briefly outline Mori's minimal model program in birational geometry.
We will discuss the relation between the naturality
problem and birational geometry in  section~4.
The main theorem was proved by a degeneration technique,
developed independently by Li-Ruan \cite{LR} and Ionel-Parker \cite{IP}.
This survey presents only Li-Ruan's approach  in  section~5.
We will discuss the possible generalization of
Mori's program in the last section. The author would like to
thank H. Clemens for explaining to him the semi-stable degenerations.
The author  benefited
from many interesting conversations with J. Koll\'ar, An-Min Li and D. Morrison on the
topics related to this article. He wishes to express his thanks to them.

\section{Quantum cohomology and its naturality}
   Due to the efforts of many people, the mathematical theory of quantum cohomology
is well-understood. We refer to \cite{R6} for a survey.
Developing the theory was rather difficult and time-consuming. Naturally, we are
looking for some applications to justify the time and energy we
spent in building the general machinery.
If we expect quantum cohomology to be as useful as cohomology,
the results so far are disappointing.

The fundamental reason that cohomology is very useful is its naturality.
Namely, a continuous map induces a  ring homomorphism on cohomology.
A fundamental problem in quantum
cohomology is the
\vskip 0.1in
\noindent
{\bf Quantum naturality problem: }{\it Define "morphism"
for symplectic manifolds so  that quantum cohomology is
natural.}
\vskip 0.1in
The author believes that the  understanding of the naturality of
quantum cohomology will be essential for the future success of
quantum cohomology theory.

   It has been known for a while that quantum cohomology is not
   natural, even with respect to  holomorphic maps. A crucial
   calculation is the quantum cohomology of projective bundles
   by Qin-Ruan \cite{QR}. The calculation clearly demonstrated
   that quantum cohomology is not natural for fibrations. This lack of naturality is
   very different from ordinary cohomology and severely limits our
   efforts to develop some nice consequences of cohomology theory,
   such as characteristic classes. The Qin-Ruan result
   shows that possible "morphisms" must be very rigid. The
   existence of these rigid morphisms will set apart quantum
   cohomology from cohomology and give it its own identity.
   The author does not know the full story at present.
   In this article, we will describe a class
   of "morphisms" of symplectic manifolds related to birational geometry.
   This class of
    "morphisms" is a certain class of surgeries called "transitions".

Let us briefly review  quantum cohomology and naturality.
We start with the definition of GW
(Gromov-Witten) invariants. Suppose that $(M, \omega)$ is
a symplectic manifold of dimension  $2n$ and $J$ is a tamed almost complex
structure. Let $\overline{\M}_A(J,g,k)$ be the moduli space of stable $J$-holomorphic genus $g$
 maps with fundamental class $A$ and $k$-marked points.
There are several methods of defining GW-invariants (\cite{FO},
\cite{LT3}, \cite{R5}, \cite{S}). One method (\cite{R5},\cite{S})
is to construct a ``virtual neighborhood'' $(U_{k}, S_{k}, E_{k})$
with the following properties: (1) $U_{k}$ is a smooth, oriented,
open orbifold whose points are maps; (2) $E_k$ is a smooth,
oriented orbifold bundle over $U_k$; (3) $S_k$ is a proper section
of $E_k$ such that $S^{-1}(0)=\overline{\M}_A(J,g,k)$. There is a
map $$ \Xi_k: U_k\rightarrow M^k $$ given by evaluating a map at
its marked points. Let $\Theta$ be a Thom form of $E_k$ supported
in a small neighborhood of the zero section. We define $$
\Psi^M_{(A,g,k)}=\int_{U_k} S^*_k \Theta\wedge \Xi^*_k \prod_i
\alpha_i\leqno(2.1) $$ for $\alpha_i\in H^*(M,\R)$. One can
eliminate divisor classes $\alpha\in H^2(M, \R)$ by the relation
$$ \Psi^M_{(A,g,k+1)}(\alpha,\alpha_1, \cdots \alpha_k)=
   \alpha(A)\Psi^M_{(A,g,k)}(\alpha_1, \cdots, \alpha_k),
\leqno(2.2)
$$
for $A\neq 0$.

The above invariants are only {\em primitive} GW-invariants.
In general, we can
also pull back cohomology classes of $\overline{\M}_{g,k}$ and insert them into the
formula (2.1).
But these more general invariants are conjectured (\cite{RT2}) to be computed
by primitive invariants.
When $A=0, g=0 k\neq 3$, we define $\Psi^M_{(A,0,k)}=0$ and
$\Psi^M_{(0,0,3)}(\alpha_1,\alpha_2, \alpha_3)=\alpha_1\cup \alpha_2\cup \alpha_3[M]$.

Choose a basis $A_1, \cdots, A_k$ of $H_2(M, \Z)$.
For $A=\sum_i a_i A_i$,
we define the formal product
$q^A=(q_{A_1})^{a_1}\cdots (q_{A_k})^{a_k}$.
For each cohomology class $w\in H^*(M)$ we  define a quantum 3-point function
$$
\Psi^M_w(\alpha_1, \alpha_2, \alpha_3)=
\sum_k\frac{1}{k!}\sum_{A} \Psi^M_{(A,0,k+3)}(\alpha_1, \alpha_2, \alpha_3,w, \cdots, w)q^A,
\leqno(2.3)
$$
where $w$ appears $k$ times.
Here, we view $\Psi^M$ as a power series in the formal variables $p_i=q^{A_i}$.
Clearly,
a homomorphism on $H_2$ will induce a change of variables $p_i$.
To define the quantum product,
we need to fix a symplectic class $[\omega]$ to define $q^A$
as an element of the Novikov ring $\Lambda_{\omega}$ (see \cite{RT1}).
Formally,  we can define quantum multiplication by the formula
$$
\alpha\times^w_Q\beta\cup \gamma[M]=\Psi^M_w(\alpha, \beta, \gamma).
\leqno(1.8)
$$
To define a homomorphism of the quantum product,
we need to match Novikov rings and hence symplectic classes.
Here, we view quantum cohomology as a theory of Gromov-Witten invariants instead of just the quantum product.
Since the GW-invariants are invariant under deformations,
the symplectic class is not a fundamental ingredient of the quantum cohomology.
In fact, the symplectic class  often
obstructs our understanding of quantum cohomology.
Clearly, the quantum 3-point function contains the same
information as the quantum product. It is  more convenient
to work directly with the quantum 3-point function.

\vskip 0.1in
\noindent
{\bf Definition 2.1: }{\it Suppose that
$$
\varphi: H_2(X, \Z)\rightarrow H_2(Y, \Z),\ H^{even}(Y, \R)\rightarrow H^{even}(X, \R)
$$
are group homomorphisms such that the maps on $H_2,H^2$
are dual to each other.
We say $\varphi$ is  {\em natural} with respect to (big) quantum cohomology if $\varphi^* \Psi^X_0=\Psi^Y_0$
($\varphi^* \Psi^X_{\varphi^*w}=\Psi^Y_{w}$)
after the change of formal variable $q^A\rightarrow q^{\varphi(A)}$.
If $\varphi$ is also an isomorphism, we say $\varphi$
induces an isomorphism on (big) quantum cohomology or $X$ and $Y$ have the same (big) quantum cohomology.}
\vskip 0.1in

 We treat two power series $F,G$ the same if $F=H+F'$, $G=H+G'$
where $G'$ is an analytic continuation of $F'$.
  For example, we can expand $\frac{1}{1-t}=\sum_{i=0}t^i$ at $t=0$
  or $\frac{1}{1-t}=\frac{1}{-t(1-t^{-1})}=-\sum_{i=0} t^{-i-1}$
  at $t=\infty$. Hence, we will treat $\sum_{i=0}t^i, \sum_{i=0} t^{-i-1}$
  as the same power series.

  When $X, Y$ are 3-folds, such a $\varphi$ is completely determined by  maps on $H_2$.
  For example, the dual map of $\varphi: H_2(X, \Z) \rightarrow H_2(Y, \Z)$ gives
  a map $H^2(Y, \R)\rightarrow H^2(X, \R)$. A map $H^4(Y, \R)\rightarrow H^4(X, \R)$
  is Poincar\'e dual to a map $H_2(Y, \R)\rightarrow H_2(X, \R)$. In the case of
  a flop, the natural map $H_2(X, \Z) \rightarrow H_2(Y, \Z)$ is an isomorphism.
  Therefore, we can take the map $H_2(Y, \R)\rightarrow H_2(X, \R)$ as its inverse.

GW-invariants are invariants under symplectic deformation. Hence, two
symplectic deformation equivalent manifolds have isomorphic big
quantum cohomology. Holomorphic symplectic varieties
are such examples. But these are trivial examples. The only known nontrivial examples
are given by the work of Li-Ruan \cite{LR} (see section 4).

   From the physical point of view, it is also natural (perhaps better) to allow
   certain {\bf Mirror Transformations}. The author does not know
   precisely what are these mirror transformations, but they appear
   naturally in mirror symmetry and in the quantum hyperplane section
   conjecture \cite{Kim}. From known examples, they  must include
   a nonlinear changes of coordinates of $H^*(X, \R), H^*(Y, \R)$ and
   a scalings of the 3-point function. Readers can find more motivation for  mirror
    transformations in the next  section.
\vskip 0.1in
\noindent
{\bf Definition 2.2: }{\it Under the assumptions of Definition 2.1, we say $\varphi$ is p-natural with respect
to (big) quantum cohomology if
we allow a mirror transformation in  Definition 2.1. In the same way, we say $\varphi$ is a p-isomorphism
of (big) quantum cohomology if
$\varphi$ is an isomorphism  and $p$-natural with respect to (big) quantum cohomology. }
\vskip 0.1in
It is obvious that naturality (isomorphism) implies p-naturality (p-isomorphism).
The author does not know any    nontrivial example of p-naturality
or p-isomorphism at present.

  Another interesting formulation is the notion of Frobenius manifold due to B. Dubrovin \cite{D}.
Here,
we require that the bilinear form $\langle\alpha, \beta\rangle=\int_X \alpha\wedge\beta$
 be
preserved under $\varphi$. This is the case for blow-downs. But the intersection form may not be preserved
under
general transitions when we have nontrivial vanishing cycles (see next section).
All three constructions (quantum products, Frobenius manifolds, quantum 3-point functions)
contain  the same quantum information. Their  difference lie in the classical information
such as the symplectic class and  bilinear form, which can be studied separately. The
author believes that Frobenius manifolds may be useful in understanding mirror
transformations.

\section{Birational geometry}

   In this section, we review Mori's minimal model program. There are several excellent reviews on
   this topic \cite{K2},\cite{K3}. We shall sketch the ideas mainly for non-algebraic geometers.
   In the process,
   we introduce the needed surgeries we need. Birational geometry is a central topic in algebraic geometry.
   The goal of birational geometry is to classify algebraic
   varieties in the same birational class. Two projective
   algebraic varieties are birational to each other iff there is a
   rational map with a rational inverse. We call this map {\em a
   birational map}. A birational map is an isomorphism between  Zariski open
   sets. Note  that a birational map is not necessarily defined everywhere.
   If it is defined everywhere, we
   call the map a contraction.
   By  definition, a contraction
   changes a lower dimensional subset only.
   Hence, we can view a contraction as a surgery.
   Intuitively, a contraction  simplifies a  variety.
   For algebraic curves, the birational equivalence is the same as the isomorphism.
    In  dimension two, any two birational algebraic
   surfaces are related by blow-up and blow-down.
   The blow-down is a contraction. If the surface
   is neither rational nor ruled, one can always perform blow-down
   until reaching a "minimal model" where blow-down can not occur.
   This is the reason that  Mori's program is often  called minimal model program.
   We want to factor a birational map as a
   sequence of contractions and find the minimal ones in the same birational class.
   In two dimensions, the minimal model is unique.
   In higher dimension, it is  much more
   difficult to carry out the minimal model program. The first difficulty is that the
   contraction is not unique.
   To overcome this difficulty, Mori
   found a beautiful correspondence between contractions and some
   combinatorial information of the algebraic manifold itself. This correspondence is
   Mori's extremal ray theory.

To define extremal rays, consider Mori's effective cone
$$
   \overline{NE}(X)=\{\sum_i a_i A_i, a_i\geq 0
    \mbox{ $A_i$ is represented by a holomorphic curve}\}
    \subset H_2(X, \R)\leqno(3.1).
$$
By definition, $\overline{NE}(X)$ is a closed cone. Let $K(X)$ be the ample cone.
One of the nice properties of projective geometry is
$$
   \overline{NE}(X)=\overline{K(X)^*}.\leqno(3.2)
$$
In symplectic geometry, we do not know if this property is true. This lack of knowledge is one
of the primary difficulties of symplectic geometry.
An edge of $\overline{NE}(X)$ is called an extremal ray of $X$.
Suppose that $L$ is an extremal ray with $K_X \cdot L<0$.
Mori showed that an extremal ray is
represented by rational curves. Each extremal ray $L$ gives a nef
class $H_L$ such that $H_L\cdot C=0$ if $[C]\in L$ and $H_L\cdot C>0$
if $[C]\not\in L$, where $C$ is a holomorphic curve. In general, {\em a class $H$ is
nef  iff $H\cdot C\geq 0$ for any holomorphic curve $C$}. Then $H_L$
defines a contraction $\phi_L: X\rightarrow Y$ by contracting every curve whose homology class is in $L$.
Moreover,  the Picard number has relation $Pic(X)=Pic(Y)+1$. $\phi_L$ is called
a primitive contraction. We abuse the notation and say that $\phi_L$ contracts $L$.
One can contract several extremal rays simultaneously.
This contraction corresponds to contracting an extremal face or other higher dimensional boundary.
 Conversely, if we have a contraction $\phi: X\rightarrow Y$ such
 that $Pic(X)=Pic(Y)+1$, then $\phi=\phi_L$ for some extremal
 ray.

Once the contraction is found, an additional difficulty arises
because $Y$ could be singular even if $X$ is smooth. Obviously, the
topology of $Y$ is simpler. However, if the singularities of $Y$
become  more and more complicated, we just trade one difficult
problem for another equally difficult problem. The next idea is
that we can  fix the class of singularities. There are two
commonly used classes of singularities: terminal singularities and canonical
singularities.
We refer readers to \cite{K2} for the precise
definitions. The terminal singularity is the minimal class of
singularities in the minimal model theory. It was shown
that in the complex dimension three, the minimal model program can
indeed be carried out. Other than for a well-understood exceptional
class of varieties (uniruled varieties), the minimal model exists.
A minimal model is, by definition, an algebraic variety $X$ such that
(i)~$X$ has terminal singularities only, and
(ii)~the canonical bundle $K_X$ is nef.

One should mention that contractions are not enough
to find the minimal model. A more difficult operation "flip" is need.
We will not give any details about flip. However, its companion "flop"
is very important to us because it resolves another difficulty of the minimal model
program. Unlike the case in dimension two, the minimal model is not unique
in  dimension three or more. However, different minimal models
are related by a sequence of flops \cite{Ka}, \cite{K1}.

We now give  a local description of a "simple flop".
Let $Y_s$ be a three-fold with a ordinary double points.
Namely, the neighborhood of a singular point
is complex analytically equivalent to a hypersurface of $\C^4$
defined by the equation
$$
    x^2+y^2+z^2+w^2=0.\leqno(3.3)
$$
The origin is the only singular point. When we blow up the origin
to obtain $Y_b$, we obtain an exceptional divisor
$E\cong \P^1\times\P^1$. Then, we contract one ruling to obtain $Y$ containing a rational
curve with normal bundle $O(-1)+O(-1)$. We can also contract
another ruling to obtain $\tilde{Y}$ containing a rational curve
with normal bundle $O(-1)+O(-1)$. Clearly, $Y, \tilde{Y}$ are the same
locally. However, they could differ globally. The process
from $Y$ to $\tilde{Y}$ is called a simple flop. For smooth
threefolds, a general flop can be deformed locally to the disjoint
union of several simple flops.

  A final difficulty is that Mori's minimal model program has only been carried out
  for complex dimension three. It is still an open problem for
  higher dimensions.

The fundamental surgeries in birational geometry are contractions,
flips, and flops. Every smooth Calabi-Yau 3-fold is a minimal model
by definition.
For minimal models, we only have flops. But flops are not enough to
classify smooth Calabi-Yau 3-folds. We need to introduce
another surgery called {\em extremal transition} or
{\em transition.}
As we mentioned previously, a contraction could lead to
a singular manifold. Sometimes, a contraction could construct a manifold with singularities  beyond the class of
terminal singularities. In fact, the flip was introduced precisely to
deal with this problem, where it will improve the singularity. The
smoothing is another well-known method to improve the singularity.  A smoothing is as follows. Consider
$$
\pi: U\rightarrow D(0,\epsilon)\subset \C\leqno(3.4)
$$
where $U$ is an analytic variety and $\pi$ is holomorphic and
of maximal rank everywhere except at zero.
Then  $X_z=\pi^{-1}(z)$ is smooth except at the central fiber $X_0$.
In this situation, we say that $X_z$ is a smoothing of $X_0$.
If $\pi$ is also of maximal rank at zero, $X_0$ is smooth and $X_z$ is deformation
equivalent to $X_0$.
Hence, they have the same quantum cohomology.
If $U$ is also K\"ahler, we say $(U, \pi)$ is a K\"ahler smoothing.
{\em A transiton is a composition of a contraction and a K\"ahler smoothing}.
We call a contraction {\em small} if the exceptional locus is of
complex codimension two or more. We call a transition {\em small}
if its corresponding contraction is small. A flip-flop is a small operation
in the sense that it only changes a subset of codimension two.
Incidentally,
flip-flop has been completely classified in dimension three.
In higher dimensions,
the classification is still an open question.
It was conjectured that
any two Calabi-Yau 3-folds can be connected to each other by a
sequence of flops or transitions and their inverses.

   We first study the change of topology regarding flops and
   transitions, which  is easy in dimension three. Suppose that $F: X\leadsto\tilde{X}$ is a
   flop, where $X, \tilde{X}$ are 3-folds. In this case, the exceptional
   locus is of complex dimension one.
   Each $A\in H_2(X, \Q)$ is represented by a pseudo-submanifold $\Sigma$.
Using PL-transversality,
   we can assume that $\Sigma$ is disjoint from the
   exceptional locus. Then $\Sigma$ can also be viewed as a pseudo-submanifold of
   $\tilde{X}$. We can also reverse this process. Therefore, $F$
   induces an isomorphism on $H_2(X, \Z)$ and hence $H^2=Hom(H_2, \Q)$.
   Using Poincar\'e duality, it induces an isomorphism on $H^4$ as
   well. The maps on $H^0, H^6$ are obvious. The first important
   theorem is

\vskip 0.1in
\noindent
{\bf Theorem 3.1 (Li-Ruan):}{\it Under the previous assumptions, $F$ induces an isomorphism on
quantum cohomology. In particular, any two three-dimensional smooth minimal models have isomorphic
quantum cohomology,  where the isomorphism is induced by flops.}
\vskip 0.1in

In dimension $\geq 4$, flop is not completely understood.
However, there are theorems of Batyrev \cite{Ba} and Wang \cite{Wa}
that any two smooth minimal models have the same Betti number. We
conjecture that
\vskip 0.1in
\noindent
{\bf Quantum Minimal Model Conjecture: }{\it Any two smooth
minimal models in any dimension have isomorphic quantum
cohomology.}
\vskip 0.1in
The Quantum Minimal Model Conjecture is true for
holomorphic symplectic varieties as well. In this case, two birational
holomorphic symplectic varieties are deformation equivalent.
Hence, there is an abstract isomorphism on big quantum cohomology.
However, it is still an interesting question to give a geometric construction of the
isomorphism.

For transitions, let $c: X\rightarrow X_s $ be the contraction
and $\pi: U\rightarrow D(0,\epsilon)$ be the K\"ahler smoothing such that $X_s$ is the central fiber.
There is a deformation retract $r: U\rightarrow X_s$. Therefore,
there is a map $\delta: H^*(X_s, \Q)\rightarrow H^*(X_z, \Q)$,
where $X_z=\pi^{-1}(z)$ is the nearby smooth fiber. In the case of
K\"ahler smoothing, the image of $\delta$ can be described in terms
of monodromy.  Readers can find
more detailed references in \cite{C}. The result is as follows. We
first blow up the singular points of $U$ along the central fiber to
obtain a smooth complex manifold $V$ whose central fibers are
a union of smooth complex submanifolds intersecting transversally with
each other. We then perform base change (pull-back by a map $z\rightarrow z^k$ on the
base) to obtain $W$. $W$ has the additional property that every
component of the central fiber has multiplicity 1 in the sense that the function
$$
\tilde{\pi}: W\rightarrow V\rightarrow U \rightarrow D(0, \epsilon)\leqno(3.5)
$$
vanishes to order one on each component. This process is called
{\em semi-stable reduction}. A consequence of semi-stable reduction is
that the monodromy $r$ of $\tilde{\pi}$ is unipotent. In this
case, we can define
$$
\log (r): H^*(X_z, \Q)\rightarrow H^*(X_z, \Q).\leqno(3.6)
$$
$\log (r)$ is a nilpotent matrix. $\ker \log (r)$ is precisely the space
of  $r$-invariant cocycles. The main result of Deligne-Schmid-Clemens \cite{C}
is that $im(\delta)^*=\ker \log (r)$. There is a natural decomposition
$$
H^*(X_z, \Q)=\ker \log (r)+im \log (r)^T.\leqno(3.7)
$$
$im \log (r)^T$ is Poincar\'e dual to the space of vanishing cycles.
Using the decomposition (3.7), we can define a map
$$
H^*(X_z, \Q)\rightarrow H^*(X_s, \Q).\leqno(3.8)
$$
Let $\phi: H^*(X_s, \Q)\rightarrow H^*(X,\Q)$ be its composition
with $c^*$. Our second theorem is that
\vskip 0.1in
\noindent
{\bf Theorem 3.2 (Li-Ruan)}{\it Suppose that $T: X\rightarrow X_z$
is a small transition in complex dimension three as described
previously. Then $\phi: H^*(X_z,\R)\rightarrow H^*(X, \R)$ is a
natural map for big quantum cohomology.}
\vskip 0.1in

The author  conjectured that
\vskip 0.1in
\noindent
{\bf Quantum Naturality Conjecture: }{\it Suppose that $T: X\rightarrow X_z$
is a small transition in any dimension and $\phi: H^*(X_z,\R)\rightarrow H^*(X, \R)$
is defined previously. Then,
\begin{description}
\item[(1)]$\Psi^{X_z}_{(A,k)}(\alpha, \cdots)=0$ for $A\neq 0$ if $\alpha$ is Poincare dual to a vanishing cycle.i.e.,
   $\alpha\in im \log(r)^T$.
   \item[(2)] Let $H\subset H^*(X_z, \R)$ be a subspace such that $\phi|_H$
   is injective. Then, $\phi|_H: H \rightarrow H^*(X,\R)$ is a
   natural map for big quantum cohomology.
\end{description}
}
\vskip 0.1in

The results above are not satisfactory because they do not say
anything about the the easiest transition, namely blow-down. The
calculation in \cite{RT1} (Example 8.6) shows that blow-down is not
natural for quantum cohomology. But this is not the end of the story.
For Calabi-Yau 3-folds, there are two other types of transitions
which are not small. These are blow-down type surgeries. By the Mirror
Surgery Conjecture \cite{LR}, there should be corresponding
operations on the Hodge structures. The operation on the
Hodge structure  is obviously natural in algebraic coordinates.
To compare with quantum cohomology, we should change
from algebraic coordinates to the flat coordinates near the large complex
structure
limit (mirror transformation). This suggests
\vskip 0.1in
\noindent
{\bf Quantum p-Naturality Conjecture: }{\it Transition induces a
p-natural map on big quantum cohomology.}
\vskip 0.1in
We hope to give shortly such an example of p-naturality in
\cite{LQR}.

For the same reason, we also have
\vskip 0.1in
\noindent
{\bf p-Quantum Minimal Model Conjecture: }{\it Any two smooth
minimal models have p-isomorphic quantum cohomology}
\vskip 0.1in
  The Quantum Minimal Model Conjecture implies the p-Quantum Minimal Model
  Conjecture; but the converse could be false.

In Mori's program, it is essential to consider varieties with
terminal singularities. It would be an important problem to work out the
singular case in complex dimension three to see if the same
phenomenon holds in the category of  varieties with terminal singularities.
Fortunately, a complete description of flops is also available in
this case \cite{K3}.

\section{Log Gromov-Witten invariants and degeneration
formula}

    In this section, we describe the techniques used to prove
Theorems 3.1, 3.2. There are two slightly different approaches by
Li-Ruan \cite{LR} and Ionel-Parker \cite{IP}. Here, we present
Li-Ruan's approach only. Readers should consult \cite{IP} for their
approach. Moreover, the author uses different terminology  from the original paper
to make it more familiar to readers who are not familiar with contact geometry.
The heart of the technique is a degeneration formula for
GW-invariants under semi-stable degeneration such that the central
fiber has only two components. The author has already given the definition
of semi-stable degeneration in the last section. In the case that
the
central fiber has only two components, the normal bundle of their intersection $Z$ has
opposite first Chern classes. During the last several years,
semi-stable degeneration reappeared in symplectic geometry under the
names ``symplectic norm sum'' and ``symplectic cutting'' and plays an
important role in some of the recent developments in symplectic
geometry. The author should point out that the recent construction
of symplectic geometry is independent from algebraic geometry. In
fact, the symplectic construction is stronger because it shows that  one can always produce a semi-stable
degeneration with prescribed central fiber in the symplectic category. Of course, the
analogous
result is far from true in the algebraic category. There are two symplectic constructions,
 norm sum and cutting, which are inverse to each other. Let me present the construction
of symplectic cutting due to E. Lerman \cite{L}.

Suppose that $H: M\rightarrow \R$ is a periodic Hamiltonian function such that the
Hamiltonian vector field $X_H$ generates a circle action. By adding a constant,
we can assume that $0$ is a regular value. Then $H^{-1}(0)$ is a smooth
submanifold preserved by the circle action. The quotient $Z=H^{-1}(0)/S^1$ is the famous
symplectic reduction. Namely, it has an induced symplectic structure.

For simplicity, we assume that $M$ has a global Hamiltonian circle action.
Once we write down the construction, we  observe that a local circle
Hamiltonian action is enough to define  symplectic cutting.

  Consider the product manifold $(M\times {\bf C}, \omega\oplus -idz\wedge
d\bar{z})$. The moment map $H-|z|^2$ generates a Hamiltonian circle action
$e^{i\theta}(x, z)=(e^{i\theta}x, e^{-i\theta}z)$. Zero is a regular value and
we have symplectic reduction
$$
\overline{M}^+=\{H=|z|^2\}/S^1.\leqno(4.1)
$$
We have the decomposition
$$
\overline{M}^+=\{H=|z|^2\}/S^1=\{H=|z|^2>0\}/S^1\cup H^{-1}(0)/S^1.\leqno(4.2)
$$
Furthermore,
$$
\phi^+: \{H>0\}\rightarrow \{H=|z|^2>0\}/S^1
$$
$$
\phi^+(x)=(x, \sqrt{H(x)}).\leqno(4.3)
$$
is a symplectomorphism.
Let
$$
M^+_b=H^{-1}(\geq 0).\leqno(4.4)
$$
Then $M^+_b$ is a manifold with boundary and there is a map
$$
M^+_b\rightarrow \overline{M}^+.\leqno(4.5)
$$
Clearly,  $\overline{M}^+$ is obtained by
collapsing the $S^1$ action on  $H^{-1}(0)$.
It is obvious that  we only need a local $S^1$ Hamiltonian action.
To obtain $\overline{M}^-$, we consider the circle action $e^{i\theta}(x,z)=(e^{i\theta}x,
e^{i\theta}z)$ with the moment map $H+|z|^2$. $\overline{M}^+, \overline{M}^-$
are called symplectic
cuttings of $M$. We define $M^-_b$ similarly. By the construction, $Z=H^{-1}(0)/S^1$ with induced
symplectic structure embedded
symplectically into $\overline{M}^{\pm}$. Moreover, its normal bundles have opposite first
Chern classes.

   Symplectic norm sum  is an inverse operation of
symplectic cutting. Gompf developed before symplectic cutting appeared.  E. Ionel  observed
that symplectic norm sum-symplectic
cutting is the same as semi-stable degeneration such that $\overline{M}^+\cup_Z \overline{M}^-$
is the central fiber.

   The main result of this section is a degeneration formula of
GW-invariants under semi-stable degeneration with the central fiber
of two components or symplectic norm sum or symplectic cutting. To
simplify the notation, we use the term symplectic cutting only.
The  first step is to introduce relative GW-invariants of a pair $(M, Z)$
where $Z$ is a smooth codimension two symplectic submanifold. We
can always choose an almost complex structure $J$ such that $Z$ is
an almost complex submanifold. We remark that
such an almost complex structure is not generic in the usual sense.
 Hence, the algebraic geometer should
view $Z$ as a smooth divisor. The process of defining relative
GW-invariants is similar to that of defining regular GW-invariants.
First, we define relative stable maps. It is helpful to recall the
definition of stable maps.
\vskip 0.1in
\noindent
{ \bf Definition 4.1 ([PW], [Ye], [KM]). }{\it A stable J-map is an equivalent class of the pair $(\Sigma,f)$.
Here,  $\Sigma$ is a nodal marked Riemann surface with arithematic
genus $g$, $k$-marked point, and $f:\Sigma \rightarrow X$ is a
continuous map whose restriction on each component of $\Sigma$ (called a component of $f$ in
short) is
$J$-holomorphic. Furthermore, it satisfies the stability
condition: if $f|_{S^2}$ is constant (called ghost bubble) for some
$S^2$-component, the $S^2$ has at least three special points
(marked points or nodal point).

If $h: \Sigma\rightarrow \Sigma'$ is a biholomorphic map,
$h$ acts on $(\Sigma, f)$ by $h\circ (\Sigma, f)=(\Sigma', f\circ
h)$. Then,  $(\Sigma, f), (\Sigma', f\circ h)$ are
equivalent or $(\Sigma,f)\sim (\Sigma',f\circ h)$.}
\vskip 0.1in

Suppose that $(\Sigma, f)$ is a $J$-stable map.
A stable map can be naturally decomposed into
connected components lying outside of $Z$ ({\em non-$Z$-factors})
or completely inside $Z$ ({\em $Z$-factors}). Both are stable
maps.  The division creates some marked points different from $x_i$.
We call these marked points {\em new marked points}.
In contrast, we call $x_i$ {\em old marked points}.
Obviously, different factors  intersect each other at new marked points.

\vskip 0.1in
\noindent
{\bf Definition 4.2: }{\it A label of $(\Sigma, f)$ consists of
(1) a division of each $Z$-factor into a set of
stable maps $\{f_{\mu}\}$ (call $Z$-subfactor)
intersecting at new marked points;
(2) an assignment of a nonzero integer $a_p$ to new marked points of
$\{f_{\mu}\}$, non-$Z$ factors with the following compatibility condition.
\begin{description}
\item[(1)] If $p$ is a new marked point of non-$Z$ factor $f_{\eta}$,
then $a_p>0$ is the order of tangency of $f_{\eta}$ with $Z$.
\item[(2)] If $p$ and $q$ are new marked points where
two components intersect, the $a_p=-a_q$.
\item[(3)] If $f_{\mu}$ and $f_{\nu}$ intersect at
$f_{\mu}(p_i)=f_{\nu}(q_i)$ for $1\leq i\leq l$,
then all the $a_{p_i}$ with $1\leq i \leq l$ have the same sign.
\end{description}
}
\vskip 0.1in

Let $N$ be the projective completion of the normal bundle
$E\rightarrow Z$, i.e., $N=P(E\oplus \C)$. Then $N$ has a zero
section $Z_0$ and an infinity section $Z_{\infty}$. We view $Z$ in
$M$ as a zero section. To define relative stable map, we assign a
nonnegative integer $t_i$ to each marked point $x_i$ such that
$\sum_i t_i=[f]Z$. Denote the tuple $(t_1, \cdots, t_k)$ by $T_k$.
\vskip 0.1in
\noindent
{\bf Definition 4.3: }{\it A log stable map
is a triple $((\Sigma, f), T, label)$ such that each $Z$-subfactor
$f_{\mu}$ can be lift to a stable map $\tilde{f}_{\mu}$ into $N$
satisfying (1) $\tilde{f}_{\mu}$ intersects $Z_0, Z_{\infty}$ at
marked points (old or new) only; (2) If $p_i$ is a new marked
point on $\tilde{f}_{\mu}$, then $\tilde{f}_{\mu}$ intersects
$Z_0$ (resp. $Z_{\infty}$) at $p_i$ with order $|a_{p_i}|$ if
$a_{p_i}>0$(resp. $a_{p_i}<0$). (3) Let $x_i$ be an old marked
point. If $x_i$ is on a non-$Z$ factor $f_{\eta}$, then $f_{\eta}$
intersects $Z$ at $x_i$ with order $t_i$. If $x_i$ is on a
$Z$-subfactor $f_{\mu}$, then $\tilde{f}_{\mu}$ intersects $Z_0$at
$x_i$ with order $t_i$.} \vskip 0.1in

It is easy to show that if the lifting $\tilde{f}_{\mu}$ exists,
it is unique up to the complex multiplication on the fiber of $N$.
Let $\overline{\M}^{M,Z}_A(g,T_k,J)$ be the moduli space of relative
stable maps with fixed $T$. Clearly, there is a map
$$\pi: \overline{\M}^{M,Z}_A(g,T_k,J)\rightarrow \overline{\M}^X_A(g,k,J).\leqno(4.6)$$

   Now let's explain the motivation of above definition. Consider
the convergence of a sequence of $J$-map $(\Sigma_n, f_n)$. Of course, $f_n$ will
converge to a stable map $(\Sigma, f)$. In general, $(\Sigma, f)$
may have some $Z$-factors. Recall that the puncture disc $D-\{0\}$
is biholomorphic to half cylinder $S^1\times [0, \infty)$. Now, we
do this fiberwisely over $Z$. We can view $M-Z$ as an almost
complex manifold with an infinite long cylinder end. We call it
{\em cylindric model}. Now, we reconsider the convergence of $(\Sigma_n,
f_n)$ in the cylindric model. The creation of a $Z$-subfactor $f_{\mu}$ corresponds to disappearance of
part of $im(f_n)$ into the infinity. We can use the translation to
rescale back missing part of $im(f_n)$. In the limit, we obtain a
stable map $\tilde{f}_{\mu}$ into $\tilde{Z}\times \R$, where
$\tilde{Z}$ is the circle bundle consisting the united vectors of $E$.
$N$ is just the closure of $\tilde{Z}\times
\R$. One can further show that $\tilde{f}_{\mu}$ indeed is a stable
map into $N$. Therefore, we obtain a lifting of $f_{\mu}$. The
label is used to specify the lifting.

Suppose that $t_i=0$ for $i\leq l$ and $t_i>0$ for $i>l$.
We have evaluation maps
$$\Xi^M_i: \overline{\B}^{M, Z}_A(g,T_k,J)\rightarrow M \leqno(4.7)$$
for $i\leq l$, and
$$\Xi^M_j: \overline{\B}^{X, Z}_A(g,T_k,J)\rightarrow Z \leqno(4.8)$$
for any $j>l$. Let $\alpha_i\in H^*(M, \R), \beta_j \in H^*(Z,
\R)$.

Roughly, the log GW-invariants are defined as $$
\Psi^{(M,Z)}_{(A,g,T_k)}(K; \alpha_1, \cdots, \alpha_l;
\beta_{l+1}, \cdots, \beta_k)=
\int_{\overline{\M}^{M,Z}_A(g,T_k,J)}\chi^*_{g,k}K\wedge
\Xi_{g,k}^*\prod_i \alpha_i\wedge P^*\prod_j \beta_j. \leqno(4.10)
$$ To be precise, the virtual techniques developed by Fukaya-Ono
\cite{FO}, Li-Tian\cite{LT3}, Ruan\cite{R5}, Siebert \cite{S}
apply to this case. For example, we can  construct a virtual
neighborhood (7.1 \cite{LR}) of $\overline{\M}^{M,Z}_A(g,T_k, J)$.
Then we integrate the integrand (1.17) over the the normalization
of  virtual neighborhood.

Clearly, there is a map $$ R:
\overline{\M}^{M,Z}_A(g,T_k,J)\rightarrow
\overline{\M}^M_A(g,k,J).\leqno(4.11) $$ However, it may not be
surjective even if $T_k=(1,\cdots, 1)$. Because of this, we remark
that when $T_k=\{1,\cdots, 1\}$, $\Psi^{(M,Z)}$ is usually
different from the ordinary or absolute GW-invariant. So
$\Psi^{(M,Z)}$ is not a ``generalized'' GW-invariant.

\vskip 0.1in

\noindent
{\bf Theorem 4.3 (Li-Ruan)(Theorem 7.6 \cite{LR}): }

\vskip 0.1in

\noindent {\it (i). $\Psi^{(M,Z)}_{(A,g,T_k)}(K; \alpha_1,...,
\alpha_l; \beta_{l+1},..., \beta_{k})$ is well-defined,
multilinear and skew-symmetric.

\vskip 0.1in

\noindent
(ii). $\Psi^{(M,Z)}_{(A,g,T_k)}(K; \alpha_1,...,\alpha_l; \beta_{l+1},...,
\beta_{k})$ is independent of the choice of forms $K, \alpha_i, \beta_j$
representing the cohomology classes $K, [\beta_j], [\alpha_i]$,  and independent of the choice
of virtual neighborhoods.

\vskip 0.1in

\noindent
(iii).
$\Psi^{(M,Z)}_{(A,g,T_k)}(K; \alpha_1,\ldots, \alpha_l; \beta_{l+1},\ldots, \beta_{k})$
is independent of the choice of almost complex
structures on $M$ where $Z$ is an almost complex submanifold, and hence is an invariant of $(M,Z)$.}

\vskip 0.1in If $K=1$, we will drop $K$ from the formula. This
invariant is called {\em a primitive log GW-invariant}. In this
article, we will give a degeneration formula for primitive
invariants and comment on how to modify it for non-primitive
invariants.

Now we explain the degeneration formula of GW-invariants under
symplectic cutting. First of all, symplectic cutting defines a map
$$
\pi: M\rightarrow \overline{M}^+\cup_Z \overline{M}^-,\leqno(4.12)
$$
where the right hand side is the central fiber of the
degeneration. By  construction, $\omega^+|_Z=\omega^-|_Z$.
Hence, the pair $(\omega^+, \omega^-)$ defines a cohomology class of $\overline{M}^+\cup_Z
\overline{M}^-$, denoted by $\omega^+\cup_Z\omega^-$. It is easy
to observe that
$$
\pi^*(\omega^+\cup_Z \omega^-)=\omega.\leqno(4.13)
$$
 $\pi$ induces a map
$$
\pi_*: H_2(M, \Z)\rightarrow H_2(\overline{M}^+\cup_Z\overline{M}^-, \Z).\leqno(4.14)
$$
Let $B\in \ker (\pi_*)$. By (4.13), $\omega(B)=0$. Define $[A]=A+\ker(\pi_*)$
and
$$
\Psi^M_{([A],g,k)}=\sum_{B\in [A]} \Psi_{(B,g,k)}.\leqno(4.15)
$$
For any $B, B'\in [A]$, $\omega(B)=\omega(B')$. By the Gromov
compactness theorem, there are only finitely many such $B$
represented by $J$-holomorphic stable maps. Therefore, the
summation in (4.15) is finite. Moreover, the cohomology class $\alpha_i$ is not arbitrary. Namely, let $\alpha^{\pm}_i\in H^*(\overline{M}^{\pm}, \R)$
such that $\alpha^+_i|_Z=\alpha^-_i|_Z$. It defines a class $\alpha^+_i\cup_Z\alpha^-_i\in H^*(\overline{M}^+\cup_Z\overline{M}^-, \R)$.
We choose $\alpha_i=\pi^*(\alpha^+\cup_Z \alpha^-)$. The degeneration formula is a big
summation
$$
\Psi^M_{([A],g,k)}(\alpha_1, \cdots, \alpha_k)=\sum_{indices} terms.\leqno(4.16)
$$
Next, we give a procedure to write down the terms on the right
hand side.

Step (1). We first write down a graph representing the topological
type of a degenerate Riemann surface in $\overline{M}^+\cup_Z\overline{M}^-$
with the following properties:
(i)~Each component is completely inside
either $\overline{M}^+$ or $\overline{M}^-$;
(ii)~the components only intersect
each other along $Z$.
(iii)~No two components where both in $\overline{M}^+$ or in $\overline{M}^-$
intersect each other;
(iv) The arithmetic genus is $g$ and the
number of marked points is $k$;
(v) The total homology class is
$\pi_*([A])$.
(vi). Each intersection point carries a positive
integer representing the order of tangency.
We use $C$ to denote such a graph.
$C$ will be used as the index of summation and (4.16)
can be written as
$$
\Psi^M_{([A],g,k)}(\alpha_1, \cdots, \alpha_k)=\sum_C \Psi_C.\leqno(4.17)
$$

Step (2). Suppose that $C$ has components $C_1, \cdots C_s$ and let $(A_{C_i}, g_{C_i})$
be the homology class and genus of  the $C_i$-component. Then,
$$
\Psi_C= r_C\prod_i \Psi^{(\overline{M}^?, Z)}_{(A_{C_i}, g_{C_i},T_i)}(variables),\leqno(4.18)
$$
where $\overline{M}^?$ is the one of $\overline{M}^{\pm}$ in which $C_i$
lies  and $r_C$ is the product of certain numbers  to be determined.
$T_i$ is given by the order of tangency. Here, the original marked
points have zero order of tangency.

Step(3). The following two steps determine the variables of  each term in (4.18)
and $r_C$. If a marked point $x_i$ appears in some component
$C_t$, then $\alpha^{\pm}_i$ should be in the variable of
$\Psi^{(\overline{M}^?, Z)}_{(A_{C_t}, g_{C_t},T_t)}(variables)$. $\pm$
depends on whether  $C_t$ lies in $\overline{M}^+$ or $\overline{M}^-$.

Step(4). $\Psi^{(\overline{M}^?, Z)}_{(A_{C_t}, g_{C_t})}(variables)$
has other variables as well associated to intersection points.
Suppose that $y$ is an intersection point of $C_i$ and $C_j$
components with order of tangency $k_y$. Let $\beta_a$ be a basis of
$H^*(Z, \R)$ and $\eta^{ab}=\int_Z \beta_a\beta_b$. Let $(\eta_{ab})$ be its inverse,
which can be
thought of as the intersection matrix of the Poincar\'e dual of $\beta_a,
\beta_b$. Then $y$ contributes a term $\beta_a$ in $\Psi_{(A_{C_i},
g_{C_i})}$, $\beta_b$ in $\Psi_{(A_{C_j}, g_{C_j})}$ and $k_y\eta_{ab}$
in $r_C$.

 The four steps above will completely determine the formula
 (4.16).
\vskip 0.1in
\noindent
{\bf Theorem 4.1 (Li-Ruan(\cite{LR}, Theorem 7.9,7.10))}{\it
There is a degeneration formula of GW-invariants under symplectic
cutting described by Step (1)-Step (4).}
\vskip 0.1in

An important comment is that different $\alpha^{\pm}_i$ may define the
same $\alpha_i$. Then we have different ways to express $\Psi^M(\cdots,\alpha_i, \cdots)$.
This is very important in the application of the degeneration formula (4.16). For example,
The Pioncar\'e dual of a point can be chosen to
have
its support completely inside $\overline{M}^+$ or
$\overline{M}^-$.

To derive a degeneration formula with class $K$, let us consider its
Poincar\'e dual $K^*$. Then we have to specify a degeneration $K^*_{\infty}$ of
$K^*$ and look at what kind of graph $C$ could appear in the
degeneration of $K^*$. Namely, we look for those $C$ for which  we obtain an element of $K^*_{\infty}$
after we contract all the unstable components of $C$. For example, suppose that $K^*=\{\Sigma\}$. We can
require that $\Sigma$ stay in the interior of $\overline{M}_{g,k}$.
This requirement will force $C$ to have one component the same as $\Sigma$ in one of $\overline{M}^{\pm}$
and  other components being  unstable rational components.

Using Theorem 4.1, we can  prove Theorem 3.1,3.2. Consider a
simple flop $X\leadsto \tilde{X}$. They have the same blow-up
$X_b$. Take a trivial family of $X$ and blow up the central fiber
along the the rational curve $L$. This is a semi-stable
degeneration whose central fiber is a union of $X_b$ and $Y_L$, where $Y_L$ is the projective completion  of
normal bundle $O(-1)+O(-1)$ of $L$ intersecting with $X_b$ along the projectivization of
the normal bundle $P(O(-1)+O(-1))\cong \P^1\times \P^1$. The crucial
information for the proof of Theorem 3.1 is that $C_1(Y_L)=3Z_{\infty}$. An index
calculation shows that the only nonzero term in (4.16) is that the $C$ lies
completely in $X_b$ or $Y_L$. We repeat the same argument for $\tilde{X}$
to express the invariant of $\tilde{X}$ in terms of that of $X_b,
Y_L$. We then obtain a formula for the change of GW-invariants
under a simple flop. To conclude the isomorphism of quantum
cohomology, we need the amazing fact that  the contribution of
multiple cover maps to $L$  cancels the error from the change of
the
classical triple product under a simple flop.

To prove Theorem 3.2, we consider the case when $X_s$ has only a
ordinary double point. We then do the semi-stable reduction to
obtain a semi-stable degeneration. In this case, the central fiber has
again only two components $X_b, Y_s$(a fact I learnt from H.
Clemens), where $X_b$ is the blow-up along $L$  and $Y_s$ is
a quadratic 3-fold intersecting with $X_b$ along a hyperplane $H$.
This is completely consistent with the
symplectic cutting description in section 2 of \cite{LR}! Note
that the first Chern class of a quadric 3-fold is $3H$! This is
precisely what we need in the case of flop. Using our degeneration
formula, we conclude that the only nonzero term in (4.16) is the case  that $C$
stays in either $X_b, Y_s$. But there is no holomophic curve in $Y_s$
not intersecting $H$. Hence, it must be in $X_b$. Therefore, we
get an expression (Theorem 8.1,\cite{LR}) of the invariant of $X_z$ in terms of the relative
invariant of $X_b$ and hence $X$ by the previous argument. The rest of the proof is
just sorting out the formula for quantum 3-fold points.

   For more general flops or small transitions, an argument of
P. Wilson shows that one can reduce it to the above case by an almost
complex deformation. Hence, our previous argument  works in
all  cases. If the reader is looking for an algebraic proof of
our theorem, it seems that one can generalize our degenerations
formula to semi-stable degeneration with an arbitrary number of
components. Then one can use this more general degeneration
formula to bypass the almost complex deformation. It is clear
that such a general degeneration formula should also contain
relative GW-invariants, but the combinatorics will be much more
complicated.

\section{Symplectic minimal model program}
   After going over  technical mathematics in the last section, it is
time to have some fun and  make some wild speculations.
Mori's minimal model program can be viewed as  a surgery theory of contractions
and flip-flops. As we mentioned in  section three, transitions
play  a crucial role in the classification of Calabi-Yau 3-folds.
The author hopes that he has convinced the reader that they also play  a crucial
role in quantum cohomology. One essence of symplectic geometry is the flexibility to deform
complex structure. Therefore, transition is also a natural operation in symplectic geometry.
It is natural that we speculate that
there should be a minimal model theory with transitions as
fundamental surgeries. We call this proposed theory the
``symplectic  minimal model program''. The holomorphic version of symplectic minimal model program
 is no longer in the category of  birational geometry.
Instead, it addresses the important problem
of connecting different moduli spaces of complex manifolds.  An interesting question is: what are the
minimal models in the holomorphic symplectic minimal model program?
Recall that in Mori  theory a minimal model has nef $K_X$.
We speculate that a minimal model  in the moduli minimal model program should be either
$Pic(X)=1$ or $K_X$ is ample. Like the case of minimal models program,
there are exceptional cases. It is of great interest
to study these exceptional cases as well.
Clearly, transition improves singularities.
The author does  not have any feeling for what kind of singularities
should be allowed in the holomorphic minimal model program.

 In my view, the guiding problem of symplectic geometry
should be the classification of  symplectic
manifolds.  Symplectic minimal model program should also be viewed as a
classification scheme of symplectic manifolds.
      In
topology, it is  rare that we can  label manifolds.
An essential step is to establish some fundamental surgeries
and classify a class of manifolds under such surgeries.
These fundamental surgeries
should simplify the manifolds. In
dimension two, the fundamental surgery is connected sum. In
dimension three, the fundamental surgeries are connected sums over
spheres and tori. In dimension 4, people are still struggling to
understand the fundamental surgeries.
A natural question is ``what are fundamental surgeries of symplectic manifolds?''
I believe that transition is
one of the fundamental surgeries of symplectic manifolds.
In fact, any  surgeries which are natural
with respect to quantum cohomology deserve  our attention, if we
believe that quantum cohomology is a fundamental invariant of
symplectic manifolds.

   Symplectic minimal model program also provide a scheme to
   attack Arnold conjectures. Recall that the Arnold conjecture
   for degenerate hamiltonian is still an open problem. The
   obstruction is precisely the existence of holomorphic curves.
   The symplectic minimal model program is designed to kill all
   the
   holomorphic curves. Therefore, Arnold conjectures should hold
   for minimal model. The rest is to do transition in Hamiltonian
   invarant fashion.

This line of thought is very appealing
because transition has a beautiful interpretation in term of
classical symplectic geometry.

   In a neighborhood of a singularity, the boundary has a natural
contact structure. It is known classically that there are two
ways to "fill" a contact manifold. A contact manifold could
bound a resolution of a singularity or a neighborhood of  the zero
section of a cotangent bundle. In the case of a singular
point, transition is a local duality which interchanges these two fillings.
In the general case, we need to consider a fiber-wise version of the above
construction.

   Therefore, it is natural to consider a "symplectic minimal model program"
   using transition as the fundamental surgery. As the author showed in \cite{R3}, one can generalize the Mori cone $NE(X)$
to symplectic manifolds. Transition simplifies a symplectic manifold $X$ in
the sense that it simplifies the Mori cone $NE(X)$.
One important ingredient in Mori's
program is the interpolation between the Mori cone and the ample
cone,
which are related by $NE(X)=\overline{K(X)^*}$. In symplectic
geometry, we no longer have such a relation. So we encounter
severe difficulties at  the first step of our symplectic
minimal model program.
It is probably a long shot to establish such a program.
But I have no doubt that much interesting
mathematics will come out of our investigation. We end our
discussion with following question: What is the minimal model in
the symplectic minimal model program?

\end{document}